\documentclass[12pt]{amsart}

\usepackage{amssymb,amsmath}
\usepackage{amsthm}
\usepackage{tikz}
\usepackage[matrix,arrow,tips]{xy}
\usepackage{color}
\usepackage{fullpage}
\usepackage{booktabs}
\usepackage{enumerate}

\usetikzlibrary{patterns}
\usetikzlibrary{calc}

\newcommand{\eps}{\epsilon}

\newcommand{\bbC}{\mathbb{C}}

\newcommand{\bbZ}{\mathbb{Z}}

\newcommand{\bbM}{\mathbb{M}}
\newcommand{\bbN}{\mathbb{N}}
\newcommand{\bbR}{\mathbb{R}}

\newcommand{\calD}{\mathcal{D}}

\newcommand{\calM}{\mathcal{M}}
\newcommand{\calN}{\mathcal{N}}

\newcommand{\calP}{\mathcal{P}}
\newcommand{\calQ}{\mathcal{Q}}

\newcommand{\calW}{\mathcal{W}}

\newcommand{\del}{\partial}

\newcommand{\End}{\mathrm{End}}

\newcommand{\Hom}{\mathrm{Hom}}

\newcommand{\diag}{\mathrm{diag}}

\newcommand{\GL}{\mathrm{GL}}

\newcommand{\Spin}{\mathrm{Spin}}
\newcommand{\SU}{\mathrm{SU}}

\newcommand{\rFs}[5]{\,_{#1}F_{#2} \left( \genfrac{.}{.}{0pt}{}{#3}{#4}
\ ;#5 \right)}

\renewenvironment{proof}{\noindent{\scshape Proof.}}{\qed}

\theoremstyle{plain}
\newtheorem{theorem}{Theorem}[section]
\newtheorem{lemma}[theorem]{Lemma}
\newtheorem{proposition}[theorem]{Proposition}
\newtheorem{corollary}[theorem]{Corollary}

\newtheorem{definition}[theorem]{Definition}

\theoremstyle{definition}
\newtheorem{example}[theorem]{Example}
\newtheorem{remark}[theorem]{Remark}

\title{Non-symmetric Jacobi polynomials of type $BC_{1}$ as vector-valued polynomials Part 1:\\ spherical functions}

\date{\today}

\author[van Horssen]{M. van Horssen}
\author[van Pruijssen]{M. van Pruijssen}

\address[van Pruijssen]{Radboud University, IMAPP-Mathematics, Heyendaalseweg 135,
6525 AJ NIJMEGEN, the Netherlands}
\email{m.vanpruijssen@math.ru.nl}

\address[van Horssen]{Radboud University, IMAPP-Mathematics, Heyendaalseweg 135,
6525 AJ NIJMEGEN, the Netherlands}
\email{max.vanhorssen@ru.nl}

\subjclass[2010]{33C52,33C45,33E30}
\keywords{non-symmetric Jacobi polynomials, spherical functions, matrix-valued orthogonal polynomials}

\begin{document}


\begin{abstract}
We study non-symmetric Jacobi polynomials of type $BC_{1}$ by means of vector-valued and matrix-valued orthogonal polynomials.
The interpretation as matrix-valued orthogonal polynomials yields a new expression of the non-symmetric Jacobi polynomials of type $BC_1$ in terms of the symmetric Jacobi polynomials of type $BC_{1}$. In this interpretation, the Cherednik operator, that has the non-symmetric Jacobi polynomials as eigenfunctions, corresponds to two shift operators for the symmetric Jacobi polynomials of type $BC_{1}$.

We show that the non-symmetric Jacobi polynomials of type $BC_{1}$ with so-called geometric root multiplicities, interpreted as vector-valued polynomials, can be identified with spherical functions on the sphere $S^{2m+1}=\Spin(2m+2)/\Spin(2m+1)$ associated with the fundamental spin-representation of $\Spin(2m+1)$. 
The Cherednik operator corresponds to the Dirac operator for the spinors on $S^{2m+1}$ in this interpretation.
\end{abstract}

\maketitle


\section{Introduction}\label{S:Intro}
The purpose of this paper is to study the non-symmetric Jacobi polynomials of type $BC_{1}$ by means of vector-valued and matrix-valued orthogonal polynomials, and to identify these polynomials with spherical functions for compact symmetric pairs.

Both symmetric and non-symmetric Jacobi polynomials have been studied in the context of a general root system $R$ with multiplicity function $k$ in \cite{MR1313912, MR1353018}, where they are defined using the Gram-Schmidt process with respect to a weight function that depends on $R$ and $k$. The non-symmetric Jacobi polynomials are eigenfunctions of the Cherednik operators, which are differential-reflection operators. The non-symmetric Jacobi polynomials can also be interpreted as vector-valued Laurent polynomials that are invariant for an action of the Weyl group of $R$. The eigenfunction property of the non-symmetric Jacobi polynomials is reflected by the corresponding vector-valued Laurent polynomials being  solutions of the KZ-equations \cite[\S3]{MR1353018}. The $q$-analogs of the (non-)symmetric Jacobi polynomials are the (non-)symmetric Macdonald polynomials \cite{MR1976581}, they have similar properties and for a suitable $q \to 1$ limit the corresponding Jacobi polynomials are recovered.

The (non-)symmetric Askey-Wilson polynomials are a special case of the (non-)symmetric Macdonald polynomials in one variable, which have the (non-)symmetric Jacobi polynomials of type $BC_{1}$ as their $q \to 1$ limit. These polynomials have been studied by Macdonald \cite[Ch.6]{MR1976581}, Stokman and Noumi \cite{MR2085854}, and more recently by Koornwinder and Bouzeffour \cite{MR3413409}. In \cite{MR3413409} the non-symmetric Askey-Wilson polynomials are expressed in terms of the symmetric Askey-Wilson polynomials and they are viewed as $\bbC^{2}$-valued polynomials. Taking the $q \to 1$ limit the authors show similar statements for the non-symmetric Jacobi polynomials. In the vector-valued interpretation the Cherednik operator corresponds to a matrix-valued differential operator that encodes the well-known forward and backward shift operators for the symmetric Jacobi polynomials \cite[(7.13)]{MR3413409}.

The key observation to describe the non-symmetric Jacobi polynomials as $\bbC^{2}$-valued polynomials is that the space $\bbC[z^{\pm1}]$, for which the non-symmetric Jacobi polynomials constitute a basis, is a free module of rank two over the algebra $\bbC[z+z^{-1}]$, for which the symmetric Jacobi polynomials constitute a basis. By choosing different generators we obtain a different family of $\bbC^{2}$-valued polynomials, whose entries can still be expressed in terms of the symmetric Jacobi polynomials. Our choice of generators leads to a new expression of the non-symmetric Jacobi polynomials in terms of the symmetric Jacobi polynomials. In this setting, the matrix-valued differential operator corresponding to the Cherednik operator encodes the contiguity shift operators for the symmetric Jacobi polynomials.

The vector-valued polynomials from \cite{MR3413409} cannot be easily combined into a family of matrix-valued orthogonal polynomials, because putting together two vector-valued polynomials of the same degree as the columns of a matrix gives a matrix-valued polynomial whose leading coefficient is not invertible. In contrast, the vector-valued polynomials that we consider allow us to produce a family of matrix-valued orthogonal polynomials. For the root multiplicity $k=(0,\nu)$ these matrix-valued orthogonal polynomials specialize to the matrix-valued Gegenbauer polynomials of size $2\times 2$ from \cite{MR3735699}.

In \cite{MR3006173,MR3085105} the columns of the matrix-valued Gegenbauer polynomials with parameter $\nu = 1$ are identified with the spherical functions for the compact symmetric pair $(\SU(2)\times\SU(2),$ $\diag(\SU(2)))$ of type $\tau_{1}$, where $\tau_{1}:\SU(2)\to\GL(\bbC^{2})$ is the standard representation. Note that $\SU(2)\cong\Spin(3)$ and $\SU(2)\times\SU(2)\cong\Spin(4)$. The fundamental spin-representation of $\Spin(2m+1)$ is an irreducible representation $\tau_{m}:\Spin(2m+1)\to\GL(V_{\tau_{m}})$ of highest weight $\omega_{m}$, where we use the standard choices of roots and weights as in \cite{MR0240238}. In this way, the triple $(\SU(2)\times\SU(2),\diag(\SU(2)),\tau_{1})$ fits into the family
$$(\Spin(2m+2),\Spin(2m+1),\tau_{m}),\quad m\ge1.$$
The compact symmetric pair $(\Spin(2m+2),\Spin(2m+1))$, with $m\ge1$, is of rank one, which implies that the spherical functions of type $\tau_{m}$ can be described with vector-valued Laurent polynomials in a single variable \cite{PvP}. Since the module generated by the spherical functions of type $\tau_{m}$ is of rank two over the algebra of zonal spherical functions, the spherical functions take values in $\bbC^{2}$, see \cite{PvP}. The compact symmetric pair gives rise to a root system, called the restricted root system, and a root multiplicity $k$. These root multiplicities encode the dimensions of the restricted root spaces and are called geometric root multiplicities. We show that the non-symmetric Jacobi polynomials of type $BC_{1}$ for geometric root multiplicities, interpreted as $\bbC^{2}$-valued Laurent polynomials, can be identified with spherical function for the pairs $(\Spin(2m+2),\Spin(2m+1))$ of type $\tau_{m}$. This provides a group-theoretic interpretation of the matrix-valued Gegenbauer polynomials with parameter $\nu\in\bbZ_{\ge0}$.

The spherical functions of type $\tau_{m}$ are in bijection with the irreducible representations of a commutative subquotient of the universal enveloping algebra of the complexified Lie algebra of $\Spin(2m+2)$, see \cite[Thm.1.4.5]{MR0954385}. This algebra acts by differential operators on the space spanned by the spherical functions, and upon identifying the spherical functions with vector-valued Laurent polynomials, these operators can be determined explicitly by so-called radial part calculations \cite{MR683007}. In fact, the spherical functions of type $\tau_{m}$ are characterized as being simultaneous eigenfunctions of the operators arising from the action of this algebra. However, for the radial part calculations explicit expressions of the corresponding elements in the universal enveloping algebra are needed, which are in general not easily obtained.
In the case that we consider, it is known that the algebra of interest is generated by a single element \cite{MR1245809}, which turns out to be the Dirac operator for spinors on the quotient $S^{2m+1}=\Spin(2m+2)/\Spin(2m+1)$, see \cite{MR1814669}. Moreover, the radial part of the Dirac operator has been calculated by  Camporesi and Pedon in \cite[Prop.5.2]{MR1814669}.

The spaces of vector-valued Laurent polynomials, where the non-symmetric Jacobi polynomials of type $BC_{1}$ and spherical functions of type $\tau_{m}$ reside, are related by the multiplication with a non-constant function. Conjugating the radial part of the Dirac operator with this non-constant function, we recover the matrix-valued differential operator corresponding to the Cherednik operator. The non-symmetric Jacobi polynomials of type $BC_{1}$ are characterized, up to normalization, as eigenfunctions of the Cherednik operator and their degree. Since this is also the case for the spherical functions of type $\tau_{m}$ with the Dirac operator, we are able to establish the desired relationship between the functions when taking their normalizations into account.

It is well-known that the symmetric Jacobi polynomials for general root systems can be identified with zonal spherical functions for compact symmetric pairs \cite[Ch.5]{MR1313912}. More recently, the intermediate Jacobi polynomials have been introduced in \cite{MvP-VVOP} as Laurent polynomials that are invariant for the action of a parabolic subgroup of the Weyl group. If the parabolic subgroup is the Weyl group itself or the trivial group then the symmetric and non-symmetric Jacobi polynomials are recovered. A natural question about intermediate Jacobi polynomials is whether they also have an interpretation as spherical functions for a compact symmetric pair $(U, K)$ of type $\tau$ and an irreducible $K$-representation $\tau$. In \cite{MvP-VVOP} it is shown that this is the case for the root system of type $A_{2}$ and the parabolic subgroup generated by a simple reflection. In this paper we show that this
is also the case for the intermediate Jacobi polynomials of type $BC_{1}$, for which the only intermediate Jacobi polynomials are the symmetric and non-symmetric Jacobi polynomials.

This paper is organized as follows. In Section \ref{S:IJP} the expressions of the non-symmetric Jacobi polynomials of type $BC_{1}$ in terms of the symmetric Jacobi polynomials of type $BC_{1}$ are derived. The identification of the non-symmetric Jacobi polynomials of type $BC_{1}$ with geometric root multiplicities as spherical functions for the pair $(\Spin(2m+2),\Spin(2m+1))$ of type $\tau_{m}$ is established in Section \ref{S:SF}.


\section{$\bbC^{2}$-valued Jacobi polynomials for the root system $BC_{1}$}\label{S:IJP}

Let $R=\{\pm\eps,\pm2\eps\}\subset\bbR^{*}$ be the root system of type $BC_{1}$, where $\bbR$ has the standard inner product (multiplication) and $\epsilon\in\bbR^{*}$ is defined by $\eps(1)=1$. The Weyl group $W$ of $R$ is the reflection group with two elements $\bbZ_{2}$. A function $R\to \bbC$ that is constant on $\bbZ_{2}$-orbits is called a root multiplicity and is determined by its values on the long and on the short roots. We represent a root multiplicity by the tuple $k=(k_{1},k_{2})$ where $k_{1}$ is the value on the short roots and $k_{2}$ on the long roots.
The weight lattice of $R$ is $P=\bbZ\eps$ and the group algebra $\bbC[P]$ is isomorphic to $\bbC[z^{\pm1}]$ by sending $e^{\eps}$ to $z$. Let
$$\delta_{k}(z)=\left(1-\frac{z+z^{-1}}{2}\right)^{k_{1}}\left(1-\frac{z^{2}+z^{-2}}{2}\right)^{k_{2}}$$
and on $\bbC[z^{\pm1}]$ the sesquilinear pairing
$$(p,q)_{k}=\int_{S^{1}}\overline{p(z)}q(z)\delta_{k}(z)\frac{dz}{iz},\quad p,q\in\bbC[z^{\pm1}]$$
which is an inner product if $k_{1},k_{2}\ge0$.
The weight lattice $P$ has a partial ordering \cite[\S3]{MvP-VVOP}, which in this one-dimensional case yields a total ordering on the monomials, $1<z<z^{-1}<z^{2}<z^{-2}<z^{3}<\cdots$. Application of the Gram-Schmidt process on this totally ordered basis of $\bbC[z^{\pm1}]$ results in an orthogonal basis $E(n,k)$ with $n\in\bbZ$ of $\bbC[z^{\pm1}]$ whose constituents have the defining properties
\begin{itemize}
\item $E(n,k)=z^{n}+$ lower order terms,
\item $(E(n,k),z^{\ell})_{k}=0$ for all monomials $z^{\ell}<z^{n}$.
\end{itemize}
The action of the Weyl group $W = \bbZ_{2} = <s>$ on $\bbC[z^{\pm1}]$ is given by $(s\cdot p)(z)=p(z^{-1})$. The Laurent polynomials $E(n,k)$ are called the non-symmetric Jacobi polynomials and they are eigenfunctions of the differential-reflection operator
$$D_{k}=z\partial_{z}+k_{1}\frac{1}{1-z^{-1}}(1-s)+2k_{2}\frac{1}{1-z^{-2}}(1-s)-\rho(k),$$
where $\rho(k)=\frac{1}{2}(k_{1}+2k_{2})$. More precisely
\begin{equation}\label{eq:BC1_eigenvalue}
\begin{aligned}
D_{k}(E(n,k))&=(n+\rho(k))E(n,k), & n>0,\\
D_{k}(E(n,k))&=(n-\rho(k))E(n,k), & n\le0.
\end{aligned}
\end{equation}
For $k_{1},k_{2}\ge0$ the eigenvalues are pairwise different and in this case \eqref{eq:BC1_eigenvalue} determines the polynomials $E(n,k)$ up to scaling. The action of the Weyl group and that of $D_{k}$ on $\bbC[z^{\pm1}]$ generate an algebra in $\End(\bbC[z^{\pm1}])$ that is isomorphic to the graded Hecke algebra of type $BC_{1}$ with root multiplicity $k$ \cite[Cor.2.9]{MR1353018}.
For later reference we record the following result.

\begin{lemma}\label{lem:HighestTermsEBC1}
Let $n\in\bbZ_{\ge0}$ and $k_{1},k_{2}\ge0$. We have
$$E(n+1,k)=z^{n+1}+\frac{k_{1}}{1+2n+2\rho(k)}z^{-n}+\mbox{lower order terms}.$$
\end{lemma} 
\begin{proof}
The condition on $k$ ensures that $1+2n+2\rho(k)\ne0$.
We have $E(n+1,k)=z^{n+1}+c_n(k)z^{-n}+\mbox{lower order terms}$ for some constant $c_n(k)$ and
$$D_{k}E(n+1,k)=(n+1+\rho(k))z^{n+1}+(-c_n(k)(n+\rho(k))+k_{1})z^{-n}+\mbox{lower order terms}$$
by direct calculation. Comparing with \eqref{eq:BC1_eigenvalue} yields $(n+1+\rho(k))c_n(k)=(-c_n(k)(n+\rho(k))+k_{1})$ from which the result follows.
\end{proof}


\subsection{Invariant $\bbC^{2}$-valued Laurent polynomials}

The Weyl group $\bbZ_{2}$ also acts on $\bbC^{2}$ by swapping the standard basis vectors $e_{1}$ and $e_{2}$. In turn $\bbZ_{2}$ acts diagonally on the space of $\bbC^{2}$-valued Laurent polynomials $\bbC[z^{\pm1}]\otimes\bbC^{2}$. The map 
\begin{equation}\label{GammaBC1}
\Gamma:\bbC[z^{\pm1}]\to\left(\bbC[z^{\pm1}]\otimes\bbC^{2}\right)^{\bbZ_{2}}, \quad\Gamma(p)=(p,s\cdot p)^T,
\end{equation}
is an isomorphism.
The space $\bbC[z^{\pm1}]\otimes\bbC^{2}$ is equipped with the sesquilinear pairing
$$(P,Q)_{k}=\frac{1}{2}\int_{S^{1}}\left(\overline{P_{1}(z)}Q_{1}(z)+\overline{P_{2}(z)}Q_{2}(z)\right)\delta_{k}(z)\frac{dz}{iz},$$
where $P,Q\in\bbC[z^{\pm1}]\otimes\bbC^{2}$ are $\bbC^{2}$-valued polynomials $P(z)=(P_{1}(z),P_{2}(z))^{T}$ and $Q(z)=(Q_{1}(z),Q_{2}(z))^{T}$.
It is an inner product if $k_{1},k_{2}\ge0$.
The map $\Gamma$ respects the sesquilinear forms. Let
\begin{equation}\label{def:VVJacobi}
P(n,k)=\Gamma(E(n,k)),\quad n\in\bbZ.
\end{equation}
Then $(P(n,k)\mid n\in\bbZ)$ is an orthogonal basis of $\left(\bbC[z^{\pm1}]\otimes\bbC^{2}\right)^{\bbZ_{2}}$.

The polynomials $P(n,k)$ with $n\in\bbZ$ are the $\bbC^{2}$-valued Jacobi polynomials with $BC_{1}$-symmetry. They are eigenfunctions of the differential operator
\begin{equation}\label{BC1:FirstD}
\Gamma_{*}(D_{k})=
\begin{pmatrix}
z\partial_{z}&0\\
0&-z\partial_{z}
\end{pmatrix}+
\begin{pmatrix}
\frac{k_{1}}{1-z^{-1}} + \frac{2k_{2}}{1-z^{-2}}-\rho(k) & -\frac{k_{1}}{1-z^{-1}} - \frac{2k_{2}}{1-z^{-2}} \\
-\frac{k_{1}}{1-z} - \frac{2k_{2}}{1-z^{2}} & \frac{k_{1}}{1-z} + \frac{2k_{2}}{1-z^{2}}-\rho(k)
\end{pmatrix}
\end{equation}
and characterized as such for $k_{1},k_{2}\ge0$.
Using the identity $\frac{1}{1-z}-1=-\frac{1}{1-z^{-1}}$, we obtain the expression
\begin{equation}\label{eqn:BC1GammaD}\Gamma_{*}(D_{k})=
\begin{pmatrix}
z\partial_{z}&0\\
0&-z\partial_{z}
\end{pmatrix}+
\left(\frac{k_{1}}{1-z^{-1}} + \frac{2k_{2}}{1-z^{-2}}\right)\begin{pmatrix}
 1 & -1 \\
1 & -1
\end{pmatrix}+
\rho(k)\begin{pmatrix}
-1&0\\
-2&1
\end{pmatrix},\end{equation}
which corresponds to the KZ-connection \cite[Def.3.1]{MR1353018}.
\begin{example}
For $k_{1},k_{2}\ge0$, the polynomials $P(0,k)$ and $P(1,k)$ are given by
\begin{equation}\nonumber
P(0,k)(z)=\begin{pmatrix}1\\1\end{pmatrix},\quad P(1,k)(z)=\begin{pmatrix} z+\frac{k_{1}}{1+2\rho(k)}\\z^{-1}+\frac{k_{1}}{1+2\rho(k)}
\end{pmatrix}.
\end{equation}
We observe that these polynomials are indeed eigenfunctions of the differential operator \eqref{eqn:BC1GammaD} with the eigenvalues $-\rho(k)$ and $1+\rho(k)$, respectively, in accord with \eqref{eq:BC1_eigenvalue}.
\end{example}


\subsection{$\bbC^{2}$-valued orthogonal polynomials}

Let $x=\frac{1}{2}(z+z^{-1})\in\bbC[z^{\pm1}]^{\bbZ_{2}}$. Note that $\bbC[z^{\pm1}]^{\bbZ_{2}} = \bbC[x]$. The $\bbC[x]$-module $\bbC[z^{\pm1}]$ is free of rank two with generators $1$ and $z$. Indeed, $\bbC[x]\cap\bbC[x]z=\{0\}$ and with induction and the formula 
$$z^n=z^{n-1}(z+z^{-1})-z^{n-2},\quad n\ge2,$$
it follows that $z^{n}\in\bbC[x]\oplus\bbC[x]z$ for all $n\in\bbZ_{\ge0}$. Since $z^{-1}=z+z^{-1}-z\in\bbC[x]\oplus\bbC[x]z$, the module $\bbC[x]\oplus\bbC[x]z$ is $\bbZ_{2}$-stable, from which we conclude
$$\bbC[z^{\pm1}]=\bbC[x]\oplus\bbC[x]z.$$
This is an isomorphism of $\bbC[x]$-modules, i.e.~for $p\in\bbC[z^{\pm1}]$ there are unique polynomials $f_{1},f_{2}\in\bbC[x]$ such that
$$p(z)=f_{1}((z+z^{-1})/2)+zf_{2}((z+z^{-1})/2).$$
We define the $\bbC[x]$-module isomorphism 
\begin{equation*}\label{CxModuleIso}
\Upsilon:\bbC[z^{\pm1}]\to\bbC[x]\otimes\bbC^{2}, \quad p=f_{1}+zf_{2}\mapsto(f_{1},f_{2})^{T}.
\end{equation*}
The space $\left(\bbC[z^{\pm1}]\otimes\bbC^{2}\right)^{\bbZ_{2}}$ is also a $\bbC[x]$-module by component-wise multiplication
$$(fP)(z)=(f((z+z^{-1})/2)P_{1}(z),f((z+z^{-1})/2)P_{2}(z))^T,\quad f\in\bbC[x],\quad P\in\left(\bbC[z^{\pm1}]\otimes\bbC^{2}\right)^{\bbZ_{2}},$$
which makes the map $\Gamma$ from \eqref{GammaBC1} a $\bbC[x]$-module isomorphism.
The $\bbC[x]$-module isomorphism
\begin{equation}\label{MapToPolyBC1}
\Upsilon\circ\Gamma^{-1}:\left(\bbC[z^{\pm1}]\otimes\bbC^{2}\right)^{\bbZ_{2}}\to\bbC[x]\otimes\bbC^{2}
\end{equation}
is given by multiplication with the inverse of the matrix-valued function
$$\Psi(z)=\begin{pmatrix}1&z\\1&z^{-1}\end{pmatrix}.$$
Let $\calW$ be the matrix-valued function 
$$\calW(x)=\begin{pmatrix}2&2x\\2x&2\end{pmatrix}.$$
Note that $\calW((z+z^{-1})/2)=\Psi(z)^{*}\Psi(z)$ on the unit circle $S^{1}\subset\bbC$, so $\calW(x)$ is a matrix-weight on the interval $[-1,1]$.
Recall \cite[Ex.1.3.2]{MR1313912} that 
$\delta_{k}(z)\frac{dz}{iz}=w_{k}(x)dx$ with
$$w_{k}(x)=2^{k_{1}+2k_{2}}(1-x)^{k_{1}+k_{2}-\frac{1}{2}}(1+x)^{k_{2}-\frac{1}{2}}.$$
The space $\bbC[x]\otimes\bbC^{2}$ is equipped with the sesquilinear pairing
\begin{equation}\label{mv-inner product}
(\calP,\calQ)_{k}=\frac{1}{2}\int_{-1}^{1}\calP(x)^{*}\calW(x)\calQ(x)w_{k}(x)dx
\end{equation}
that is an inner product for $k_{1},k_{2}\ge0$. The map \eqref{MapToPolyBC1} respects the sesquilinear forms. It will be convenient to look at the matrix-valued polynomials that are orthogonal with respect to the weight $\calW(x)w_{k}(x)$.


\subsection{$\bbM_{2}$-valued orthogonal polynomials}

Let $\bbM_{2}$ be the algebra of $2\times 2$-matrices with complex entries. The space of $\bbM_{2}$-valued polynomials is denoted by $\bbM_{2}[x]$. For two polynomials $\calP,\calQ\in\bbM_{2}[x]$ the formula \eqref{mv-inner product} defines an $\bbM_{2}$-valued inner product, see \cite[\S1]{MR3298987}.

\begin{definition}\label{def:VVJacobiX}
For $n\in\bbZ$ let $\calP(n,k)\in\bbC[x]\otimes\bbC^{2}$ be defined by
$$\calP(n,k)((z+z^{-1})/2)=\Psi(z)^{-1}P(n,k)(z)$$ for $n\in\bbZ$, where $P(n,k)$ is defined in \eqref{def:VVJacobi}.

For $N\in\bbZ_{\ge0}$ let $\calM(N,k)$ be the $\bbM_{2}$-valued polynomial whose first column is $\calP(-N,k)$ and second column $\calP(N+1,k)$. 
\end{definition}

Note that $\calP(n,k)$ is the image of $P(n,k)$ under the map $\Upsilon\circ\Gamma^{-1}$ and hence a polynomial in $x=(z+z^{-1})/2$.
Likewise, $(\calM(N,k)\mid N\in\bbZ_{\ge0})$ is a family of $\bbM_{2}$-valued orthogonal polynomials for the matrix-weight $\calW w_{k}$. We proceed to calculate the polynomials $\calM(N,k)$. 

\begin{lemma}
The polynomial $\calM(N,k)$ is of degree $N$ and its leading coefficient is
\begin{equation}\label{BC1-CN}
C_{N}(k)=2^{N}\begin{pmatrix}
1& \frac{k_{1}}{1+2N+2\rho(k)}\\0&1\end{pmatrix}.
\end{equation}
\end{lemma}

\begin{proof}
In view of Definition \ref{def:VVJacobiX} we have 
\begin{equation}\label{BC1-M}
\calM(N,k)((z+z^{-1})/2)=\begin{pmatrix}1&z\\1&z^{-1}\end{pmatrix}^{-1} \begin{pmatrix}
E(-N,k)(z)& E(N+1,k)(z)\\  E(-N,k)(z^{-1})  & E(N+1,k)(z^{-1})
\end{pmatrix}.
\end{equation}
Lemma \ref{lem:HighestTermsEBC1} implies $E(N+1,k)(z)=z^{N+1}+c_{N}(k)z^{-N}+$ lower order terms. Likewise $E(-N,k)(z)=z^{-N}+d_{N}(k)z^{N}+$lower order terms, for some coefficient $d_{N}(k)$. Using
$$\Psi(z)^{-1}=\frac{1}{z-z^{-1}}\begin{pmatrix}-z^{-1}&z\\1&-1\end{pmatrix},$$ 
a straightforward calculation shows that $\calM(N,k)(x)=C_{N}(k)x^{N}+$lower order terms.
\end{proof}

Note that we can also conjugate the differential operator \eqref{eqn:BC1GammaD} by $\Psi(z)$ and change the variable to $x=(z+z^{-1})/2$ to obtain the differential operator
$$
\calD_{k}=\begin{pmatrix}
-x&-1\\1&x
\end{pmatrix}\del_{x}+\begin{pmatrix}
-\rho(k)&k_{1}\\0&1+\rho(k)
\end{pmatrix},$$
where $\partial_{x}$ is component-wise differentiation with respect to $x$. The polynomials $\calM(N,k)$ are eigenfunctions of $\calD_{k}$ in the following sense,
\begin{equation*}\label{BC1DIffOpM}
\calD_{k}(\calM(N,k))=\calM(N,k)\Lambda(N,k),\quad \Lambda(N,k)=\begin{pmatrix}
-N-\rho(k)&0\\0&N+1+\rho(k)
\end{pmatrix}.
\end{equation*}


\subsection{Diagonalization of the matrix-weight}

The matrix-weight $\calW$ diagonalizes by a constant matrix,
$$U\calW U^{*}=2\begin{pmatrix}1-x&0\\0&1+x\end{pmatrix},\quad U=\frac{1}{\sqrt{2}}\begin{pmatrix}1&-1\\1&1
\end{pmatrix}.$$
It follows that $\{U\calM(N,k)\mid N\in\bbZ_{\ge0}\}$ is a family of $\bbM_{2}$-valued orthogonal polynomials for the weight
\begin{equation}\label{DiagJacobiBC1}
2^{k_{1}+2k_{2}+1}
\begin{pmatrix}
(1-x)^{k_{1}+k_{2}+\frac{1}{2}}(1+x)^{k_{2}-\frac{1}{2}}&0\\
0&(1-x)^{k_{1}+k_{2}-\frac{1}{2}}(1+x)^{k_{2}+\frac{1}{2}}
\end{pmatrix}.
\end{equation}
By multiplication on the right with the leading coefficient we obtain
\begin{equation}\label{BC1monicM}
U\calM(N,k)C_{N}(k)^{-1}U^{-1}=x^{N}+\mbox{lower order terms},\quad N\in\bbZ_{\ge0},
\end{equation}
which is a monic family of matrix-valued orthogonal polynomials for the weight \eqref{DiagJacobiBC1}.
At the same time, we note that the matrix-weight \eqref{DiagJacobiBC1} is diagonal with Jacobi weights on the diagonal entries. Define
\begin{equation}\label{BC1alphabeta}
\alpha=k_{1}+k_{2}-\frac{1}{2},\quad\beta=k_{2}-\frac{1}{2}
\end{equation}
and let
$$P^{(\alpha,\beta)}_{n}(x)=\frac{(\alpha+1)_{n}}{n!}\rFs{2}{1}{-n,n+\alpha+\beta+1}{\alpha+1}{\frac{1-x}{2}}$$
be the classical Jacobi polynomial of degree $n$ with parameters $\alpha,\beta$. 
Recall that $P^{(\alpha,\beta)}_{n}((z+z^{-1})/2)$ is a multiple of
$E(n,k)(z)+E(n,k)(z^{-1})$, the symmetric Jacobi polynomial of degree $n$ for the root system $BC_{1}$ \cite[Thm.2.12]{MR1353018}.
The matrix-weight \eqref{DiagJacobiBC1} becomes
\begin{equation}\label{DiagJacobiBC1alphabeta}
2^{\alpha+\beta+2}
\begin{pmatrix}
(1-x)^{\alpha+1}(1+x)^{\beta}&0\\
0&(1-x)^{\alpha}(1+x)^{\beta+1}
\end{pmatrix}
\end{equation}
and
$$\calN(N,(\alpha,\beta))=\begin{pmatrix}P^{(\alpha+1,\beta)}_{N}&0\\0&P^{(\alpha,\beta+1)}_{N}\end{pmatrix}
,\quad N\in\bbZ_{\ge0},$$
is a family of matrix-valued orthogonal polynomials for \eqref{DiagJacobiBC1alphabeta}. Note that 
\begin{equation}\label{BC1-N}
\calN(N,(\alpha,\beta))(x)=\left(\frac{1}{2}\right)^{N}\frac{(N+\alpha+\beta+2)_{N}}{N!}x^{N}+\mbox{lower order terms}.
\end{equation}

\begin{proposition}\label{prop:relations}
For $N\in\bbZ_{\ge0}$, $k_{1},k_{2}\ge0$ and the convention \eqref{BC1alphabeta} we have
\begin{multline*}\calM(N,k)(x)
=\frac{2^{2N-1}N!}{(N+\alpha+\beta+2)_{N}}\cdot \\
\begin{pmatrix}
P_{N}^{(\alpha+1,\beta)}(x) + P_{N}^{(\alpha,\beta+1)}(x) & -P_{N}^{(\alpha+1,\beta)}(x) + P_{N}^{(\alpha,\beta+1)}(x)\\
-P_{N}^{(\alpha+1,\beta)}(x) + P_{N}^{(\alpha,\beta+1)}(x) & P_{N}^{(\alpha+1,\beta)}(x) + P_{N}^{(\alpha,\beta+1)}(x)
\end{pmatrix}
\begin{pmatrix}
1& -\frac{\alpha-\beta}{2N+\alpha+\beta+2}\\0&1\end{pmatrix}.
\end{multline*}
\end{proposition}

\begin{proof}
The monic family of orthogonal $\bbM_{2}$-valued polynomials for a given matrix-weight is unique, so an expression of $\calM(N,k)$ is obtained by comparing \eqref{BC1monicM} with \eqref{BC1-N}. The result follows from \eqref{BC1-M}.
\end{proof}

As a corollary of Proposition \ref{prop:relations} we obtain expressions for the non-symmetric polynomials in terms of the symmetric Jacobi polynomials. 
\begin{corollary}\label{cor:relations}
Let $N\in\bbZ_{\ge0}$, $k_{1},k_{2}\ge0$ and denote $c_{N}(k)=\frac{k_{1}}{1+2N+2\rho(k)}$ for the $(1,2)$-coefficient in \eqref{BC1-CN}.
Then, with $x=(z+z^{-1})/2$ and the convention \eqref{BC1alphabeta}, we have
\begin{multline}\label{E-N}
E(-N,k)(z)=\\
\frac{2^{2N-1}N!}{(N+\alpha+\beta+2)_{N}}\left(P_{N}^{(\alpha+1,\beta)}(x)+P_{N}^{(\alpha,\beta+1)}(x)-z(P_{N}^{(\alpha+1,\beta)}(x)-P_{N}^{(\alpha,\beta+1)}(x))\right)
\end{multline}
and
\begin{multline}\label{EN1}
E(N+1,k)(z)=\\
\frac{z2^{2N-1}N!}{(N+\alpha+\beta+2)_{N}}\left(P_{N}^{(\alpha+1,\beta)}(x)+P_{N}^{(\alpha,\beta+1)}(x)
-c_{N}(k)(P_{N}^{(\alpha+1,\beta)}(x)-P_{N}^{(\alpha,\beta+1)}(x))\right)\\
+\frac{2^{2N-1}N!}{(N+\alpha+\beta+2)_{N}}\left(c_{N}(k)(P_{N}^{(\alpha+1,\beta)}(x)+P_{N}^{(\alpha,\beta+1)}(x))-P_{N}^{(\alpha+1,\beta)}(x)+P_{N}^{(\alpha,\beta+1)}(x)\right).
\end{multline}
\end{corollary}

\begin{remark}\label{Rem:KB}
An expression of the non-symmetric Jacobi polynomials in terms of the symmetric Jacobi polynomials has also been derived in \cite[\S7]{MR3413409} by different means, and in the Appendix we show that these expressions are essentially the same as those in Corollary \ref{cor:relations}. 
\end{remark}

The differential operator
\begin{equation*}\label{BC1:frakD}
\mathfrak{D}_{(\alpha,\beta)}=U\calD_{k}U^{-1}=\begin{pmatrix}
0&-(1+x)\\1-x&0\end{pmatrix}\del_{x}+\frac{1}{2}\begin{pmatrix}
\beta-\alpha+1&-2-2\beta\\
-2-2\alpha&\alpha-\beta+1
\end{pmatrix}
\end{equation*}
has $\calN(N,(\alpha,\beta))$ as an eigenfunction by construction,
\begin{equation}\label{BC1-otherShift}
\mathfrak{D}_{(\alpha,\beta)}\calN(N,(\alpha,\beta))=\calN(N,(\alpha,\beta))\mathfrak{L}(N,(\alpha,\beta)),
\end{equation}
where $\mathfrak{L}(N,(\alpha,\beta))=UC_{N}(k)\Lambda(N,k)C_{N}(k)^{-1}U^{-1}$.
We have
$$\mathfrak{L}(N,(\alpha,\beta))=\frac{1}{2}\begin{pmatrix}
\beta-\alpha+1&-2\beta-2N-2\\
-2\alpha-2N-2&\alpha-\beta+1
\end{pmatrix}.$$
This implies the following well-known result.
\begin{corollary}\label{cor:shift}
The Jacobi polynomials satisfy
\begin{align*}
((x+1)\partial_{x}+\beta+1)P^{(\alpha,\beta+1)}_{N}(x)&=(\beta+1+N)P^{(\alpha+1,\beta)}_{N}(x),\\
((x-1)\partial_{x}+\alpha+1)P^{(\alpha+1,\beta)}_{N}(x)&=(\alpha+1+N)P^{(\alpha,\beta+1)}_{N}(x),
\end{align*}
which encode the contiguity shift operators \cite[(3.3.4,5)]{MR1313912}.
\end{corollary}

\begin{proof}
Note that \eqref{BC1-otherShift} implies
$$\begin{pmatrix}-(1+x)\del_{x}P^{(\alpha,\beta+1)}_{N}\\(1-x)\partial_{x}P^{(\alpha+1,\beta)}_{N}
\end{pmatrix}=\begin{pmatrix}
(\beta+1)P^{(\alpha,\beta+1)}_{N}-(\beta+1+N)P^{(\alpha+1,\beta)}_{N}\\
(\alpha+1)P^{(\alpha+1,\beta)}_{N}-(\alpha+1+N)P^{(\alpha,\beta+1)}_{N}
\end{pmatrix},
$$
from which the result follows.
\end{proof}

The relation \eqref{BC1-otherShift} and its connection to the contiguity shift operators have been found for the root multiplicity $k=(0,0)$ in \cite[\S8.2.1]{MR3006173}.

We find triples $(\bbM_{2}[x],\calW w_{k},\calD_{k})$ consisting of a space matrix-valued polynomials $\bbM_{2}[x]$, a matrix-weight $\calW w_{k}$ and a differential operator $\calD_{k}$ that has a family of $\bbM_{2}$-valued orthogonal polynomials as simultaneous eigenfunctions, which are in turn orthogonal with respect to the $\bbM_{2}$-valued inner product obtained by integrating against the matrix-weight $\calW w_{k}$ over the interval $[-1,1]$. 

We see that although the matrix-weight diagonalizes, the differential operator does not. The corresponding orthogonal polynomials are diagonal, but the diagonal entries are coupled by the differential operator. It is interesting to find shift operators in disguise in this way. 


\section{Spherical functions}\label{S:SF}

Consider the family of compact symmetric pairs
$$(U_{m},K_{m})=(\Spin(2m+2),\Spin(2m+1)),\quad m=1,2,3,\ldots$$
with $U_{m}/K_{m}=S^{2m+1}$ the corresponding symmetric spaces.
Note that these pairs are strong Gelfand pairs, i.e.~each irreducible representation of $K_{m}$ induces multiplicity free to $U_m$, see e.g.~\cite{MR2598998}.
Following the structure theory for symmetric pairs, let $A_{m}\subset U_{m}$ be the (one-dimensional) torus for which, among other things, $U_{m}=K_{m}A_{m}K_{m}$. Let $M_{m}=Z_{K_{m}}(A_{m})$ be the centralizer of $A_{m}$ in $K_{m}$. The irreducible representations of $U_{m}$, $K_{m}$, and $M_{m}$ are denoted by $\pi$, $\tau$, and $\sigma$, respectively. There are two fundamental spin-representations of $U_{m}$ and $M_{m}=\Spin(2m)$ that we denote by $\pi_{m}^{\pm}$ and $\sigma_{m}^{\pm}$, respectively. In contrast, $K_{m}$ has only one fundamental spin-representation that is denoted by $\tau_{m}$. The restriction of $\tau_{m}$ to $M_{m}$ splits into the direct sum of $\sigma_{m}^{+}$ and $\sigma_{m}^{-}$. It follows that the spherical functions for the pair $(U_{m},K_{m})$ of type $\tau_{m}$, upon restriction to $A_{m}$, take values in a two-dimensional space.

To describe the irreducible representations of $U_{m}$ that contain $\tau_{m}$ upon restriction to $K_{m}$ we identify for the irreducible representations with their highest weights, where we follow \cite{MR0240238} for the standard choices of tori and roots and weights to make the identification:
\begin{itemize}
\item $\varpi_{1},\ldots,\varpi_{m+1}$ denote the fundamental weights of $\Spin(2m+2)$,
\item $\omega_{1},\ldots,\omega_{m}$ denote the fundamental weights of $\Spin(2m+1)$,
\item $\eta_{1},\ldots,\eta_{m}$ denote the fundamental weights of $\Spin(2m)$.
\end{itemize}

Given irreducible unitary representations $\tau:K_{m}\to\GL(V_{\tau})$ and $\pi:U_{m}\to\GL(V_{\pi})$ with $\dim(\Hom_{K_{m}}(V_{\pi},V_{\tau}))=1$, let
$$j:V_{\tau}\to V_{\pi},\quad p:V_{\pi}\to V_{\tau}$$
be non-trivial unitary $K_{m}$-intertwiners for which $p\circ j=\mathrm{Id}|_{V_{\tau}}$. The spherical function of type $\tau$ associated to $\pi$ is the matrix-valued function
$$\Psi^{\tau}_{\pi}:U_{m}\to\End(V_{\tau}),\quad \Psi^{\tau}_{\pi}(u)=p\circ\pi(u)\circ j.$$
Note that $\Psi^{\tau}_{\pi}(e)=\mathrm{Id}_{V_{\tau}}$. If $\tau = \tau_0$ is the trivial representation then the spherical functions are called zonal spherical functions.

The weight lattice for $U_{m}$ is denoted by $P_{U_{m}}$ and the subset of dominant weights by $P^{+}_{U_{m}}$. If we want to specify the highest weight then we indicate this with a subscript, e.g.~$\pi_{\varpi_{1}}$ denotes the irreducible representation of $U_{m}$ of highest weight $\varpi_{1}$.

Since $(U_{m},K_{m})$ is a Gelfand pair, we have $\dim\Hom_{K_{m}}(V_{\tau_{0}},V_{\pi})\le1$ with equality precisely if $\pi=\pi_{\ell\varpi_{1}}$ with $\ell\in\bbN_{0}$, cf.~\cite{MR0528837}. These irreducible representations of $U_{m}$ are called spherical representations and their highest weights are collected in the set
$$P^{+}_{U_{m}}(0)=\{\ell\varpi_{1}\mid\ell\in\bbZ_{\ge0}\}.$$
The spin-representation $\tau_{m}=\tau_{\omega_{m}}$ is of highest weight $\omega_{m}$ and it also has the property that $\dim\Hom_{K_{m}}(V_{\tau_{m}},V_{\pi})\le1$. In this case we have equality if and only if $\pi$ is an irreducible representation whose highest weight is contained in the set
$$P^{+}_{U_{m}}(\omega_{m})=\{\varpi_{+},\varpi_{-}\}+P^{+}_{U_{m}}(0),$$
where $\varpi_{-}=\varpi_{m}$ and $\varpi_{+}=\varpi_{m+1}$ are the fundamental spin-weights of $U_{m}$.
This can be proved by the classical branching rules but also from \cite{PvP} where $m\ge 3$ corresponds to item B5 in \cite[Table 2]{PvP} and $m=1,2$ to items B11 and B1.3 in \cite[Table 2]{PvP}, respectively.

Let $\Psi_{\pi}^{\tau_{m}}$ denote the spherical function on $U_{m}$ of type $\tau_{m}$ associated to $\pi$.
The restriction of $\Psi_{\pi}^{\tau_{m}}$ to $A_{m}$ takes values in $\End_{M_{m}}(V_{\tau_{m}})$, which is isomorphic to $\bbC^{2}$ by sending the block corresponding to $\sigma_{m}^{+}$ to $e_{1}$ and the block corresponding to $\sigma_{m}^{-}$ to $e_{2}$. Using these identifications, we denote the restricted spherical function $\Psi_{\pi_{\varpi_{\pm}+\ell\varpi_{1}}}^{\tau_{m}}|_{A_{m}}$ by
\begin{equation}\label{psiplusminus}
\Psi_{\pm}(\ell)=\begin{pmatrix}\psi_{\pm}^{+}(\ell)\\ \psi_{\pm}^{-}(\ell)
\end{pmatrix},
\end{equation}
implicitly defining the functions $\psi_{\pm}^{\pm}(\ell)$ on $A_{m}$. The algebra of regular functions on $A_{m}$ is denoted by $\bbC[A_{m}]$ which makes $\Psi_{\pm}(\ell)$ and element of $\bbC[A_{m}]\otimes\bbC^{2}$.

Let $G_{m}=\Spin(2m+2,\bbC)$ and $H_{m}=\Spin(2m+1,\bbC)$ be the complexifications of $U_{m}$ and $K_{m}$. Let $A_{m,\bbC}$ be the complexification of $A_{m}$.
The spherical functions for the non-compact Cartan dual $(\Spin(2m+1,1),\Spin(2m+1))$ of type $\tau_{m}$ have been investigated in \cite{MR1814669}. Such a spherical function is determined by its restriction to $A_{m}^{nc}=A_{m,\bbC}\cap\Spin(2m+1,1)$.
The spherical functions of type $\tau_{m}$ that are associated to the principal series representations with parameter $\lambda\in\bbC\setminus\{\pm1\}$, when restricted to $A_{m}^{nc}$, are denoted by
\begin{equation}\label{phiplusminus}
\Phi_{\pm}(\lambda)=\begin{pmatrix}\phi_{\pm}^{+}(\lambda)\\ \phi_{\pm}^{-}(\lambda)
\end{pmatrix}.
\end{equation}
Note that in \cite{MR1814669} the transformation behavior of the spherical functions is inverted, and our $\Phi_{\pm}(\lambda)$ correspond to those in \cite[(5.9)]{MR1814669} with an inverted argument. Writing $\varphi_{\pm}^{\pm}(\lambda)$ for the scalar components of \cite[(5.9)]{MR1814669}, we have $\phi_{\pm}^{\pm}(\lambda)(z)=\varphi_{\pm}^{\pm}(\lambda)(z^{-1})=\varphi_{\pm}^{\mp}(\lambda)(z)$, see the comment above (5.6) in \cite{MR1814669} for the second equality.
The regular functions on $A_{m}$ correspond to the holomorphic functions on $A_{m,\bbC}$.

The spherical functions of type $\tau_{m}$ are simultaneous eigenfunctions of a commutative subquotient of the universal enveloping algebra of the complexified Lie algebra of $\Spin(2m+2)$,\linebreak see \cite[Thm.1.4.5]{MR0954385}. This algebra also acts on the space of spinors and it is known to be generated by the Dirac operator, see \cite[Thm.1]{MR1245809} or \cite[Thm.4.3(II)]{MR1814669}.

\begin{lemma}\label{CampPed}
Let $m\ge1$ be an integer. The restriction to $A_{m}$ of a spherical function for the pair $(U_{m},K_{m})$ of type $\tau_{m}$ is an eigenfuncion of the differential operator
$$R_{m}=\begin{pmatrix}
z\partial_{z}&0\\
0&-z\partial_{z}
\end{pmatrix}+m\left(2\frac{z^{2}+z^{-2}}{z^{2}-z^{-2}}\begin{pmatrix} 1&0\\0&-1\end{pmatrix}
+\frac{4}{z^{2}-z^{-2}}\begin{pmatrix}0&-1\\1&0
\end{pmatrix}\right),
$$
with eigenvalue $\pm2i\lambda(\ell)$.
\end{lemma}

\begin{proof}
The functions $\psi_{\pm}^{+}(\ell)$ and $\psi_{\pm}^{-}(\ell)$ from \eqref{psiplusminus} can be extended to $A_{m,\bbC}$ and restricted to $A_{m}^{nc}$, and then they match with the scalar components $\phi_{\pm}^{+}(\lambda(\ell))$ and $\phi_{\pm}^{-}(\lambda(\ell))$ of the spherical functions $\Phi_{\pm}(\lambda(\ell))$ for some $\lambda(\ell) \in \bbC\setminus\{\pm1\}$. The existence of $\lambda(\ell)$ follows from Casselman's subrepresentation theorem \cite[Prop.4.2.3]{MR0929683}, and the explicit $\ell$-dependence of $\lambda(\ell)$ is not necessary for the rest of the proof but is established in Corollary \ref{cor:lambdamatch} below.

The functions $\phi_{\pm}^{\pm}(\lambda(\ell))=\varphi_{\pm}^{\mp}(\lambda(\ell))$ satisfy the system of equations of \cite[Prop.5.2.II]{MR1814669}. However, this choice of coordinate on $A_{m}^{nc}$ is not suitable to describe spherical functions on $A_{m}$ because $M_{m}\cap A_{m}\cong\bbZ/2\bbZ$ is non-trivial and has to be taken into account, while $M_{m}\cap A_{m}^{nc}$ is trivial. The system of equations of \cite[Prop.5.2.II]{MR1814669} with $w=\exp(t e_{1})$ and $\frac{d}{dt}=w\partial _{w}$ reads
\begin{equation}\label{EqInW}
\begin{aligned}
\left(w\partial_{w}+m\frac{w+w^{-1}}{w-w^{-1}}\pm i\lambda(\ell)\right)\varphi_{\pm}^{+}(\lambda(\ell))(w)&=\frac{2m}{w-w^{-1}}\varphi_{\pm}^{-}(\lambda(\ell))(w),\\
\left(w\partial_{w}+m\frac{w+w^{-1}}{w-w^{-1}}\mp i\lambda(\ell)\right)\varphi_{\pm}^{-}(\lambda(\ell))(w)&=\frac{2m}{w-w^{-1}}\varphi_{\pm}^{+}(\lambda(\ell))(w),
\end{aligned}
\end{equation}
and upon changing to $z=\sqrt{w}, w\partial_{w}=\frac{1}{2}z\partial_{z}$ we see that \eqref{EqInW} is equivalent to
\begin{equation*}\label{EqInZ}
\begin{aligned}
\left(z\partial_{z}+2m\frac{z^{2}+z^{-2}}{z^{2}-z^{-2}}\mp 2i\lambda(\ell)\right)\psi_{\pm}^{+}(\ell)(z)=\frac{4m}{z^{2}-z^{-2}}\psi_{\pm}^{-}(\ell)(z),\\
\left(z\partial_{z}+2m\frac{z^{2}+z^{-2}}{z^{2}-z^{-2}}\pm 2i\lambda(\ell)\right)\psi_{\pm}^{-}(\ell)(z)=\frac{4m}{z^{2}-z^{-2}}\psi_{\pm}^{+}(\ell)(z),
\end{aligned}
\end{equation*}
using that $\psi_{\pm}^{\pm}(\ell)(z)=\varphi_{\pm}^{\mp}(\lambda(\ell))(z^2)$.
\end{proof}

We proceed to show that the spherical functions can be identified with the non-symmetric Jacobi polynomials.
The degree of the Laurent polynomials in the entries of the restricted spherical function $\Psi_{\pm}(\ell)$ is maximal in the entry corresponding to the $M_{m}$-type $\sigma_{m}^{\mp}$, i.e.~
\begin{align*}
\Psi_{+}(\ell)(z)&=a_{+}(\ell)\begin{pmatrix}0\\z^{-(2\ell+1)}\end{pmatrix}+
\mbox{lower order terms},\\
\Psi_{-}(\ell)(z)&=a_{-}(\ell)\begin{pmatrix}z^{-(2\ell+1)}\\0\end{pmatrix}+
\mbox{lower order terms},
\end{align*}
with the coefficients $a_{\pm}(\ell)\ne0$, by the proof of \cite[Lem.6.1]{MR3801483}.

\begin{lemma}
The fundamental spherical functions of type $\tau_{m}$ associated to $\varpi_{\pm}$ are
$$\Psi_{+}(0)(z)=\begin{pmatrix}z\\z^{-1}\end{pmatrix},\quad \Psi_{-}(0)(z)=\begin{pmatrix}z^{-1}\\z\end{pmatrix}.$$
In particular, the restrictions to $A_{m}$ of these spherical functions are independent of $m$.
\end{lemma}

\begin{proof}
The function $\Psi(z)=az^{-1}+bz+c$, with $a,b,c\in\bbC^{2}$, is an eigenfunction of $R_{m}$ if and only if $\Psi=\Psi_{\pm}(0)$ and the corresponding eigenvalue is $\pm(2m+1)$.
Since these eigenfunctions are of the right degree, they must be the indicated spherical functions.
\end{proof}

Let $E_{A_{m}}\subset\bbC[A_{m}]\otimes\bbC^{2}$ be the vector space spanned by the restrictions to $A_{m}$ of the spherical functions of type $\tau_{m}$. It is a free module of rank two over the ring of $K_{m}$-biinvariant functions on $U_{m}$ restricted to $A_{m}$, see \cite[Thm.8.12]{PvP}, the generators are $\Psi_{\pm}(0)$. Note that the space of $K_{m}$-biinvariant functions restricted to $A$ is $\bbC[z^{\pm2}]^{\bbZ_{2}}$.
\begin{lemma}\label{lemma: EAm=1}
The space $E_{A_{m}}$ is isomorphic to $(\bbC[z^{\pm2}]\otimes\bbC^{2})^{\bbZ_{2}}$ via the multiplication with
$$\begin{pmatrix}
z&0\\0&z^{-1}
\end{pmatrix}.
$$
\end{lemma}
\begin{proof}
The generators $\Psi_{\pm}(0)$ can be rewritten into
$$\begin{pmatrix}
z^{-1}&0\\0&z
\end{pmatrix}\begin{pmatrix}z^{2}\\z^{-2}\end{pmatrix},\quad \begin{pmatrix}
z^{-1}&0\\0&z
\end{pmatrix}\begin{pmatrix}1\\1\end{pmatrix}.$$
So the indicated map sends generators of $E_{A_{m}}$ to generators of $(\bbC[z^{\pm2}]\otimes\bbC^{2})^{\bbZ_{2}}$, both viewed as $\bbC[z^{\pm2}]^{\bbZ_{2}}$-modules.
\end{proof}

Under this isomorphism the operator $R_{m}$ is conjugated into the operator $Q_{m}$ on $(\bbC[z^{\pm2}]\otimes\bbC^{2})^{\bbZ_{2}}$ given by
$$Q_{m}=\begin{pmatrix} z&0\\0&z^{-1}\end{pmatrix}\circ R_{m} \circ \begin{pmatrix} z^{-1}&0\\0&z\end{pmatrix},$$
and more explicitly by
$$Q_{m}=\begin{pmatrix}
z\partial_{z}&0\\
0&-z\partial_{z}
\end{pmatrix}-\begin{pmatrix}1&0\\0&1\end{pmatrix}+m\left(2\frac{z^{2}+z^{-2}}{z^{2}-z^{-2}}\begin{pmatrix} 1&0\\0&-1\end{pmatrix}
+\frac{4}{z^{2}-z^{-2}}\begin{pmatrix}0&-z^{2}\\z^{-2}&0
\end{pmatrix}\right).$$

To compare with the theory of non-symmetric Jacobi polynomials we consider the restricted root system for the compact symmetric pair $(U_{m},K_{m})$ which is $\Sigma'=\{\pm\eps\}$ of type $B_{1}$ with $\eps/2$ the fundamental weight and the root multiplicity being $2m$. We view this as a subsystem of $\Sigma=\{\pm\eps/2,\pm\eps\}$ of type $BC_{1}$ with root multiplicity $2m$. The coordinate on the torus is $z=e^{\eps/2}$ as in Lemma \ref{CampPed}.

Now consider $R=2\Sigma=\{\pm\eps,\pm2\eps\}$ with root multiplicity $k=(0,m)$. The coordinate $\zeta=e^{\eps}$ that we use, is equal to $\zeta=z^{2}$. The corresponding differential operator \eqref{BC1:FirstD} reads
$$
2\Gamma_{*}(D_{k})=
\begin{pmatrix}
z\partial_{z}&0\\
0&-z\partial_{z}
\end{pmatrix}+
\begin{pmatrix}
\frac{4m}{1-z^{-4}}-2m &  -\frac{4m}{1-z^{-4}} \\
-\frac{4m}{1-z^{4}} & \frac{4m}{1-z^{4}}-2m
\end{pmatrix}=Q_{m}+I.
$$

\begin{theorem}\label{thm:spinors}
The spherical functions $\Psi_{\pm}(\ell)(z)$ for the pairs $(\Spin(2m+2),\Spin(2m+1))$ of type $\tau_{m}$ can be identified with the non-symmetric Jacobi polynomials for the root system $2\Sigma$ with root multiplicity $(0,m)$ via
\begin{equation}\label{eqn:DegreeMatch}
\begin{aligned}
\Psi_{+}(\ell)(z)&=\frac{(\ell+2m+1)_{\ell}}{2^{2\ell}(m+\frac{1}{2})_{\ell}}\begin{pmatrix}
z^{-1}&0\\0&z
\end{pmatrix}P(\ell+1,(0,m))(z^2),\quad \ell=0,1,2,\ldots,\\
\Psi_{-}(\ell)(z)&=\frac{(\ell+2m+1)_{\ell}}{2^{2\ell}(m+\frac{1}{2})_{\ell}}\begin{pmatrix}
z^{-1}&0\\0&z
\end{pmatrix}P(-\ell,(0,m))(z^2),\quad \ell=0,1,2,\ldots.
\end{aligned}
\end{equation}
\end{theorem}

\begin{proof}
The normalization follows from Corollary \ref{cor:relations} in conjunction with \cite[(4.1.1)]{MR0372517}.
For $\ell\ge0$ we have
\begin{align*}
P(\ell+1,(0,m))(z^2)&=\begin{pmatrix}0\\z^{-(2\ell+2)}\end{pmatrix}+\mbox{lower order terms},\\
P(-\ell,(0,m))(z^2)&=\begin{pmatrix}z^{-2\ell}\\0\end{pmatrix}+\mbox{lower order terms},
\end{align*}
which follows from a small calculation using Lemma \ref{lem:HighestTermsEBC1} and the fact that $k_{1}=0$. Since $\Psi_{\pm}(\ell)(1)=(1,1)^{T}$, we conclude that leading terms on both sides of the equations in \eqref{eqn:DegreeMatch} match. At the same time, the functions on both sides of the equations in \eqref{eqn:DegreeMatch} are eigenfunctions of $R_{m}$ and hence they are the same.
\end{proof}

As a corollary of Theorem \ref{thm:spinors} we determine the eigenvalues of $R_{m}$. The eigenvalues cannot be determined from the group theory, in contrast with the eigenvalues of the Casimir operator, because there is no analog of the Harish-Chandra homomorphism available in this context. However, the $\ell$-dependence of $\lambda(\ell)$ can be traced back using the eigenvalues.
\begin{corollary}\label{cor:lambdamatch}
We have
\begin{align*}
R_{m}\left(\Psi_{+}(\ell)\right)&=(2\ell+2m+1)\Psi_{+}(\ell),\quad \ell\ge0,\\
R_{m}\left(\Psi_{-}(\ell)\right)&=-(2\ell+2m+1)\Psi_{-}(\ell),\quad \ell\ge0,
\end{align*}
and the $\ell$-dependence of $\lambda(\ell)$ in Lemma \ref{lemma: EAm=1} is given by $\lambda(\ell)=-i(\ell+m+\frac{1}{2})$.
\end{corollary}
The Jacobi polynomials are special instances of the Jacobi functions $\phi^{(\alpha, \beta)}_\lambda$, i.e.~we have
\begin{equation}\label{eq:jf}
\phi_{-i(2\lambda + \alpha + \beta + 1)}^{(\alpha, \beta)}(t)= \frac{\ell!}{(\alpha + 1)_\ell} P^{(\alpha, \beta)}_\ell(\cosh(2t)),
\end{equation}
see \cite[(2.4)]{MR0774055}.
In conjunction with our expressions for the nonsymmetric Jacobi polynomials from Corollary \ref{cor:relations} we recover the expressions for the spherical functions from \cite{MR1814669}. 

\begin{corollary}\label{cor:match}
For $n = 2m + 1$ and $z^2 = e^t$ we have
\begin{align*}
\varphi_\pm^+(\lambda)(t) &= \cosh(t/2) \phi_{2\lambda}^{(n/2 - 1, n/2)}(t/2) \mp i \frac{2\lambda}{n} \sinh(t/2) \phi_{2\lambda}^{(n/2, n/2 - 1)}(t/2),\\
\varphi_\pm^-(\lambda)(t) &= \cosh(t/2) \phi_{2\lambda}^{(n/2 - 1, n/2)}(t/2) \pm i \frac{2\lambda}{n} \sinh(t/2) \phi_{2\lambda}^{(n/2, n/2 - 1)}(t/2),
\end{align*}
for all $\lambda=\lambda(\ell)$ with $\ell=0,1,2,\ldots$. These equations correspond to (5.24) and (5.25) from \cite{MR1814669}, respectively.
\end{corollary}

\begin{proof}
Writing $R^{(\alpha, \beta)}_\ell(x) = \frac{\ell!}{(\alpha + 1)_\ell} P^{(\alpha, \beta)}_\ell(x)$, \eqref{E-N} and \eqref{EN1} imply
$$z^{-1} \frac{(\ell+n)_{\ell}}{2^{2\ell}(n/2)_{\ell}}E(-\ell, (0,m))(z^2)=\frac{z + z^{-1}}{2} R_{\ell}^{(n/2-1,n/2)}(x)+\frac{2\ell+n}{n}\frac{z-z^{-1}}{2} R_{\ell}^{(n/2,n/2-1)}(x),$$
and
$$z^{-1} \frac{(\ell+n)_{\ell}}{2^{2\ell}(n/2)_{\ell}}E(\ell+1,(0,m))(z^2)=\frac{z + z^{-1}}{2} R_{\ell}^{(n/2-1,n/2)}(x)-\frac{2\ell+n}{n}\frac{z-z^{-1}}{2} R_{\ell}^{(n/2,n/2-1)}(x),$$
respectively. By Theorem \ref{thm:spinors},
$$\varphi_\pm^-(\lambda(\ell)) = \psi_\pm^+(\ell) = \frac{(\ell+n)_{\ell}}{2^{2\ell}(n/2)_{\ell}}z^{-1} E(\ell+1, (0, m)),$$
which implies the second relation in conjunction with \eqref{eq:jf}. The first relation follows from this upon interchanging $z$ and $z^{-1}$.
\end{proof}

\begin{remark}
The spherical functions of Theorem \ref{thm:spinors} with root multiplicity $k=(0,1)$ have been obtained in \cite{MR3006173,MR3085105}, which is based on the work of Koornwinder \cite{MR783984}. In \cite{MR3735699} and in \cite{vPR2023} the multiplicity $k=(0,1)$ is extended to the parameter $k=(0,\nu)$, without the interpretation of spherical functions. 
\end{remark}


\appendix\section{}

To compare the expressions from Corollary \ref{cor:relations} with \cite[\S7]{MR3413409}, see \eqref{KB75}  and \eqref{KB74} below, we introduce some notation. We still use the convention \eqref{BC1alphabeta} to relate the parameters $(k_{1},k_{2})$ and $(\alpha,\beta)$. Define
\begin{align*}
E_{N}[z;\alpha,\beta]&=\frac{2N+\alpha+\beta+1}{N+\alpha+\beta+1}E(-N,k)(z^{-1}),\\
E_{-(N+1)}[z;\alpha,\beta]&=2E(N+1,k)(z^{-1}),
\end{align*} 
for all $N\ge0$. Furthermore, let 
$$P_{N}[z;\alpha,\beta]=\frac{2^{2N}N!}{(N+\alpha+\beta+1)_{N}}P^{(\alpha,\beta)}_{N}((z+z^{-1})/2).$$
The relations from \cite[\S7]{MR3413409} are
\begin{align}
E_{N}[z;\alpha,\beta]&=P_{N}[z;\alpha,\beta]+\frac{N}{N+\alpha+\beta+1}(z-z^{-1})P_{N-1}[z,\alpha+1,\beta+1],\label{KB75}\\
E_{-(N+1)}[z;\alpha,\beta]&=P_{N+1}[z;\alpha,\beta]-(z-z^{-1})P_{N}[z,\alpha+1,\beta+1]\label{KB74}.
\end{align}

The contiguous relations for Jacobi polynomials imply
\begin{align}
P^{(\alpha+1,\beta)}_{N}(x)-P^{(\alpha,\beta+1)}_{N}(x)&=P^{(\alpha+1,\beta+1)}_{N-1}(x),\label{contig1}\\
P^{(\alpha+1,\beta)}_{N}(x)+P^{(\alpha,\beta+1)}_{N}(x)&=2P^{(\alpha,\beta)}_{N}(x)+xP^{(\alpha+1,\beta+1)}_{N-1}(x).\label{contig2}
\end{align}
which can be used to rewrite \eqref{E-N} into
\begin{multline*}
E(-N,k)(z)=\frac{2^{2N-1}N!}{(N+\alpha+\beta+2)_{N}}\left(
2P^{(\alpha,\beta)}_{N}(x)+xP^{(\alpha+1,\beta+1)}_{N-1}(x)
-z
P^{(\alpha+1,\beta+1)}_{N-1}(x)
\right)\\
=\frac{2^{2N-1}N!}{(N+\alpha+\beta+2)_{N}}\left(
2P^{(\alpha,\beta)}_{N}(x)-\frac{z-z^{-1}}{2}P^{(\alpha+1,\beta+1)}_{N-1}(x)
\right)\\
=\frac{N+\alpha+\beta+1}{2N+\alpha+\beta+1}\left(P_{N}[z;\alpha,\beta]-(z-z^{-1})\frac{N}{N+\alpha+\beta+1}P_{N-1}[z;\alpha+1,\beta+1]\right),
\end{multline*}
which implies \eqref{KB75}. For the relation \eqref{KB74} plug the relations \eqref{contig1} and \eqref{contig2} into \eqref{EN1},
\begin{multline*}
E(N+1,k)(z)=\\
\left(x+\frac{z-z^{-1}}{2}\right)\frac{2^{2N-1}N!}{(N+\alpha+\beta+2)_{N}}\left(
2P_{N}^{(\alpha,\beta)}(x)+xP_{N-1}^{(\alpha+1,\beta+1)}(x)
-c_{N}(k)
P_{N-1}^{(\alpha+1,\beta+1)}(x)
\right)\\
+\frac{2^{2N-1}N!}{(N+\alpha+\beta+2)_{N}}\left(c_{N}(k)(
2P_{N}^{(\alpha,\beta)}(x)+xP_{N-1}^{(\alpha+1,\beta+1)}(x)
)
-P_{N-1}^{(\alpha+1,\beta+1)}(x)
\right)
\end{multline*}
and write it as a sum of symmetric and anti-symmetric term for $z\leftrightarrow z^{-1}$,
\begin{multline}\label{KB74prep}
\frac{(N+\alpha+\beta+2)_{N}}{2^{2N-1}N!}E(N+1,k)(z)=\\
(2x+2c_{N}(k))P_{N}^{(\alpha,\beta)}(x)+(x^2-1)P_{N-1}^{(\alpha+1,\beta+1)}(x)\\
+(z-z^{-1})\left(
P_{N}^{(\alpha,\beta)}(x)+\frac{1}{2}(x-c_{N}(k))P_{N-1}^{(\alpha+1,\beta+1)}(x)
\right).
\end{multline}
The identity \eqref{KB74} is equivalent to \eqref{KB74prep}, which can be seen by comparison of the even and odd parts and the following result.
\begin{lemma} The following identities of Jacobi polynomials hold,
\begin{align}
\frac{4(N+1)}{2N+\alpha+\beta+2}P^{(\alpha,\beta)}_{N+1}(x)&=(2x+2c_{N}(k))P_{N}^{(\alpha,\beta)}(x)+(x^2-1)P_{N-1}^{(\alpha+1,\beta+1)}(x)\label{Eq:J1}\\
\frac{N+\alpha+\beta+2}{2N+\alpha+\beta+2}P^{(\alpha+1,\beta+1)}_{N}(x)&=P_{N}^{(\alpha,\beta)}(x)+\frac{1}{2}(x-c_{N}(k))P_{N-1}^{(\alpha+1,\beta+1)}(x).\label{Eq:J2}
\end{align} 
\end{lemma}

\begin{proof}
These identities follow from \cite[(4.5.7)]{MR0372517} together with the forward and backward shift relations \cite[(9.8.7),(9.8.8)]{MR2656096}.
\end{proof}


\textbf{Acknowledgment.} We thank Tom Koornwinder, Erik Koelink, and the anonymous referees for useful remarks on an earlier version of this manuscript. 


\bibliography{MvH-and-MvP-Non-sym-BC1-part-1-v2}{}
\bibliographystyle{plain}

\end{document}